\documentclass[12pt]{article}

\usepackage{amsmath,amssymb,amsbsy,amsfonts,amsthm,latexsym,
                    amsopn,amstext,amsxtra,euscript,amscd}

\begin{document}

\newtheorem{theorem}{Theorem}
\newtheorem{lemma}[theorem]{Lemma}
\newtheorem{claim}[theorem]{Claim}
\newtheorem{cor}[theorem]{Corollary}
\newtheorem{prop}[theorem]{Proposition}
\newtheorem{definition}{Definition}
\newtheorem{question}[theorem]{Open Question}

\def\cA{{\mathcal A}}
\def\cB{{\mathcal B}}
\def\cC{{\mathcal C}}
\def\cD{{\mathcal D}}
\def\cE{{\mathcal E}}
\def\cF{{\mathcal F}}
\def\cG{{\mathcal G}}
\def\cH{{\mathcal H}}
\def\cI{{\mathcal I}}
\def\cJ{{\mathcal J}}
\def\cK{{\mathcal K}}
\def\cL{{\mathcal L}}
\def\cM{{\mathcal M}}
\def\cN{{\mathcal N}}
\def\cO{{\mathcal O}}
\def\cP{{\mathcal P}}
\def\cQ{{\mathcal Q}}
\def\cR{{\mathcal R}}
\def\cS{{\mathcal S}}
\def\cT{{\mathcal T}}
\def\cU{{\mathcal U}}
\def\cV{{\mathcal V}}
\def\cW{{\mathcal W}}
\def\cX{{\mathcal X}}
\def\cY{{\mathcal Y}}
\def\cZ{{\mathcal Z}}

\def\A{{\mathbb A}}
\def\B{{\mathbb B}}
\def\C{{\mathbb C}}
\def\D{{\mathbb D}}
\def\E{{\mathbb E}}
\def\F{{\mathbb F}}
\def\G{{\mathbb G}}
\def\I{{\mathbb I}}
\def\J{{\mathbb J}}
\def\K{{\mathbb K}}
\def\L{{\mathbb L}}
\def\M{{\mathbb M}}
\def\N{{\mathbb N}}
\def\O{{\mathbb O}}
\def\P{{\mathbb P}}
\def\Q{{\mathbb Q}}
\def\R{{\mathbb R}}
\def\S{{\mathbb S}}
\def\T{{\mathbb T}}
\def\U{{\mathbb U}}
\def\V{{\mathbb V}}
\def\W{{\mathbb W}}
\def\X{{\mathbb X}}
\def\Y{{\mathbb Y}}
\def\Z{{\mathbb Z}}

\def\ep{{\mathbf{e}}_p}

\def\scr{\scriptstyle}
\def\\{\cr}
\def\({\left(}
\def\){\right)}
\def\[{\left[}
\def\]{\right]}
\def\<{\langle}
\def\>{\rangle}
\def\fl#1{\left\lfloor#1\right\rfloor}
\def\rf#1{\left\lceil#1\right\rceil}
\def\le{\leqslant}
\def\ge{\geqslant}
\def\eps{\varepsilon}
\def\mand{\qquad\mbox{and}\qquad}

\def\vec#1{\mathbf{#1}}
\def\inv#1{\overline{#1}}
\def\vol#1{\mathrm{vol}\,{#1}}

\def\Ker{\operatorname{Ker}}
\def\Im{\operatorname{Im}}
\def\Re{\operatorname{Re}}
\def\deg{\operatorname{deg}}
\def\det{\operatorname{det}}
\def\tr{\operatorname{tr}}
\def\End{\operatorname{End}}
\def\Ext{\operatorname{Ext}}
\def\Hom{\operatorname{Hom}}
\def\Aut{\operatorname{Aut}}
\def\Gal{\operatorname{Gal}}
\def\Ind{\operatorname{Ind}}
\def\Centr{\operatorname{Centr}}
\def\Frob{\operatorname{Frob}}
\def\Cl{\operatorname{Cl}}
\def\id{\operatorname{id}}
\def\mod{\operatorname{mod}}
\def\dim{\operatorname{dim}}
\def\disc{\operatorname{disc}}
\def\car{\operatorname{char}}
\def\exp{\operatorname{exp}}
\def\lcm{\operatorname{lcm}}
\def\gcd{\operatorname{gcd}}
\def\e{\operatorname{e}}
\def\Idemp{\operatorname{Idemp}}
\def\GL{\operatorname{GL}}
\def\SL{\operatorname{SL}}
\def\PGL{\operatorname{PGL}}
\def\PSL{\operatorname{PSL}}
\def\li{\operatorname{li}}
\def\O{\operatorname{O}}
\def\o{\operatorname{o}}
\def\h{\operatorname{h}}
\def\log{\operatorname{log}}
\def\genus{\operatorname{genus}}
\def\tors{\operatorname{tors}}
\def\sep{\operatorname{sep}}
\def\e{\varepsilon}
\def\sp{{\operatorname{sp}}}
\def\Tr{\operatorname{Tr}}
\def\Nm{\operatorname{Nm}}
\def\Fr{\operatorname{Fr}}
\def\mymod{{\;\mbox{mod}\;}}

\newcommand{\comm}[1]{\marginpar{%
\vskip-\baselineskip 
\raggedright\footnotesize
\itshape\hrule\smallskip#1\par\smallskip\hrule}}

\def\xxx{\vskip5pt\hrule\vskip5pt}


\title{\bf Distribution of Farey Fractions in Residue
Classes and Lang--Trotter Conjectures on Average}

\author{
{\sc Alina Carmen Cojocaru} \\
{Department of Math., Stat. and Comp. Sci.}\\
{University of Illinois at Chicago} \\
{Chicago, IL, 60607, USA} \\
{cojocaru@math.uic.edu}
\\
\\
{\sc Igor E. Shparlinski} \\
{Department of Computing, Macquarie University} \\
{Sydney, NSW 2109, Australia} \\
{igor@ics.mq.edu.au}
}
\date{\today}
\pagenumbering{arabic}

\maketitle

\begin{abstract} We prove that the set of Farey fractions
of order $T$, that is, the set
$\{\alpha/\beta \in \Q\ : \ \gcd(\alpha, \beta) = 1, \ 1 \le \alpha, \beta \le T\}$,
is uniformly distributed in residue classes modulo a
prime $p$ provided $T  \ge p^{1/2 +\eps}$ for any fixed $\eps>0$.
We apply this  to obtain
upper bounds for the Lang--Trotter conjectures on Frobenius traces and
Frobenius fields
``on average'' over a one-parametric family of elliptic curves.
\end{abstract}


\section{Introduction}
\label{sec:intro}

\hspace*{0.5cm}
For a real positive  $T$, we consider the  set of Farey fractions
$$
\cF(T) = \{\alpha/\beta \in \Q\ : \ \gcd(\alpha, \beta) = 1, \ 1 \le
\alpha, \beta \le T\},
$$
for which we know that
$$
\# \cF(T) = \(\frac{6 }{\pi^2} + o(1) \)T^2.
$$
We use some results of~\cite{Shp} to
show that the elements of this set are uniformly distributed
in residue classes modulo a prime  $p$.
More precisely, we prove:

\begin{theorem}
\label{thm:R(T)}
Let $p$ be a fixed prime.
For an integer $v$, we denote by $R_{T,p}(v)$ the number
of fractions $\alpha/\beta \in \cF(T)$ with $\gcd(\beta,p)=1$
and $\alpha/\beta \equiv v \pmod p$.
Then
$$
\sum_{ 1\leq v \leq p-1  }
  \left|R_{T,p}(v) - \frac{6 }{\pi^2} \cdot \frac{T^2}{p}\right|
=
O\(T^2p^{-1} + T p^{1/2 + o(1)}\).
$$
\end{theorem}

We  apply this result to  study  Lang--Trotter conjectures
``on average'' for specializations at elements of $\cF(T)$ of
the elliptic curve
\begin{equation}
\label{eq:Family AB}
E(t) : \quad Y^2 = X^3 + A(t)X + B(t)
\end{equation}
over $\Q(t)$, where  $A(t), B(t) \in \Z[t]$.
For a general background on elliptic curves, we refer the reader
to~\cite{Silv}.
To state our results for elliptic curves, let us first recall some
standard notation.

Given an elliptic curve $E$ over $\Q$ and $a \in \Z$,
we denote by
$\Pi_{E} (a, x)$ the number of primes $p \le x$ which do not
divide the conductor $N_E$  of $E$ and such that
$a_p(E) = p + 1 -\# E_p(\F_p) = a$, where $E_p$ denotes
 the reduction $E_p$ of $E$ modulo $p$.

For a fixed imaginary quadratic field $\K$,
we denote by $\Pi_{E}(\K, x)$ the number of primes $p \le x$ which
do not divide $N_E$ and such that $a_p(E) \neq 0$ and
$\Q\left(\sqrt{a_p(E)^2 - 4 p}\right) = \K$.

Two celebrated Lang--Trotter conjectures assert that: if $a \neq 0$,
or $a=0$ and $E$ is without complex multiplication,
then
$$
\Pi_{E}(a, x) =  \(c(E,a) +o(1)\) \frac{\sqrt{x}}{\log x}
$$
for some constant $c(E,a) \ge 0$ depending only on $E$ and $a$;
if $E$ is without complex multiplication, then
$$
\Pi_{E}(\K, x) =  \(C(E,\K) +o(1)\) \frac{\sqrt{x}}{\log x}
$$
for some constant $C(E,\K) > 0$ depending only on $E$ and $\K$.

Despite a series of  interesting (conditional and unconditional)
results, these conjectures are widely open;
even the Extended Riemann Hypothesis
only allows to obtain  upper bounds on $\Pi_{E}(a, x)$
(and lower bounds in the case $a=0$)
and $\Pi_{E}(\K, x)$,
and  those are  not of the
conjectured order of magnitude; see, for
example,~\cite{CojDav,CoFoMu, El, FoMu, MurMurSar,Ser},
and also the recent surveys~\cite{Coj,MuShp}.
Therefore it makes sense to study $\Pi_{E}(a, x)$
and $\Pi_{E} (\K, x)$ on average over
some natural families of curves.
For example, Fouvry and Murty~\cite{FoMu}, and
David and Pappalardi~\cite{DavPapp1},  have considered
the average of $\Pi_{E}(a, x)$
for the family of curves $Y^2 = X^3 + uX + v$
where the integers $u$ and $v$ satisfy the inequalities
$|u| \le U$, $|v|\le V$; they have shown that if
$UV\ge x^{3/2 + \eps}$ and $\min\{U,V\} \ge x^{1/2 + \eps}$
for some fixed  positive $\eps>0$, then ``on average'' the Lang-Trotter
conjecture holds for such curves. This result has been
extended in various
directions~\cite{AkDavJur,Baier,BaSh,BBIJ,DavPapp2,Gek,James,JamYu}.
Cojocaru and Hall~\cite{CojHal} have recently considered the
one parametric family of curves of the form~\eqref{eq:Family AB}
and established an improved upper bound on the average value of
$\Pi_{E}(a, x)$ over curves of such families when
the parameter $t$ runs through the elements of $\cF(T)$
with $T$ of the same order of magnitude as $x$.

Since obtaining tight  ``individual'' estimates
is an ultimate goal, it also makes sense to
reduce the amount of  ``averaging''.
In this direction, we show that
one can obtain the bound  of~\cite[Theorem~4]{CojHal}
(established for $T \gg x$)
starting already with  $T \ge x^{3/4+\eps}$ for some fixed $\eps > 0$.
We also obtain a similar result for $\Pi_{E} (\K, x)$.
More precisely, we prove:

\begin{theorem}
\label{thm:Aver L-T}
Let $A(t), B(t) \in \Z[t]$ be fixed polynomials such that
$E(t)$ given by~\eqref{eq:Family AB}
is an elliptic curve over $\Q(t)$ with non-constant $j$-invariant,
that is,
$$
\Delta(t) = -16 (4A(t)^3 + 27B(t)^2) \neq 0
$$
and
$$
j(t)  = -\frac{6912A(t)^{3}}{4A(t)^3 + 27B(t)^2} \not\in \Q.
$$
Then for arbitrary real positive  $x$ and $T$,
\begin{enumerate}
\item
for any integer  $a \neq 0$,
$$
\sum_{\substack{\tau \in \cF(T)\\ \Delta(\tau) \ne 0}}
\Pi_{E(\tau)}(a, x)
\ll T^2 x^{3/4} + T  x^{3/2 + o(1)};
$$
\item for  $a = 0$,
$$
\sum_{\substack{\tau \in \cF(T)\\ \Delta(\tau) \ne 0}}
\Pi_{E(\tau)}(0, x)
\ll T^2 x^{2/3} + T  x^{3/2 + o(1)};
$$
\item for any imaginary quadratic field $\K$,
$$
\sum_{\substack{\tau \in \cF(T)\\ \Delta(\tau) \ne 0}}
\Pi_{E(\tau)} (\K, x)
\ll T^2 x^{2/3} + T  x^{3/2 + o(1)}.
$$
\end{enumerate}
\end{theorem}

It is easy to see that, for $T\ge x^{3/4+\varepsilon}$ for any fixed $\eps>0$,
the bound of Part~1 of Theorem~\ref{thm:Aver L-T} becomes
$$
\frac{1}{\#\cF(T)} \sum_{\substack{\tau \in \cF(T)\\ \Delta(\tau) \ne 0}}
\Pi_{E(\tau)}(a, x)
\ll x^{3/4};
$$
this is exactly the same as the bound  of~\cite[Theorem~4]{CojHal},
which, however,  had been established only for $T\gg x$.
Similarly,  for $T\ge x^{5/6 +\varepsilon}$ for any fixed $\eps>0$,
the bounds of Parts~2 and 3 become
$$
\frac{1}{\#\cF(T)} \sum_{\substack{\tau \in \cF(T)\\ \Delta(\tau) \ne 0}}
\Pi_{E(\tau)} (0, x)
\ll x^{2/3};
$$
$$
\frac{1}{\#\cF(T)} \sum_{\substack{\tau \in \cF(T)\\ \Delta(\tau) \ne 0}}
\Pi_{E(\tau)}(\K, x)
\ll x^{2/3}.
$$
We also see that Theorem~\ref{thm:Aver L-T} is nontrivial
for  $T\ge x^{1/2+\varepsilon}$.

We recall
   that the notations $U \ll V$ and  $U = O(V)$ are both equivalent to
the statement that $|U| \le c V$ holds with some constant $c> 0$,
which throughout the paper may depend on the polynomials $A(t)$ and
$B(t)$ in~\eqref{eq:Family AB}.
We also use $o(1)$ to  denote a quantity which tends to zero
as $T\to \infty$.

\bigskip

{\bf Acknowledgements.} The authors are very grateful to Chris Hall
for valuable comments.

 The first author would like to thank the Fields
Institute for an excellent working environment.
This work was supported in part  by NSF grant DMS 0636750 (for A.C. ~C.)
and by ARC grant DP0556431 (for I.~S.).


\section{Proof of Theorem~\ref{thm:R(T)}}

\hspace*{0.5cm}
First,  we note that $R_{T,p}(0) = O(T^2/p)$, which is within the total
error term of Theorem~\ref{thm:R(T)}.
Thus  it is enough
to concentrate on $R_{T, p}(v)$ with  $v =1, \ldots, p-1$.

For an integer $d$ we let
\begin{eqnarray*}
M_{W,p,d}(v) =
\#
\{(\alpha, \beta) \in \Z^2 \  : &  
 1\le \alpha,\beta \le W,\ d\mid
\gcd(\alpha,\beta),\ \gcd(p,\beta) =1,
\\
& \qquad \qquad \qquad \qquad \quad 
 \alpha/\beta \equiv v \pmod p\}.
\end{eqnarray*}
Clearly,  $M_{W,p,d}(t) = 0$ if $p|d$ and
$
M_{W,p,d}(t) = M_{W/d,p,1}(t).
$

Now  let  $\mu(d)$ denote the M\"obius function.
Using  the inclusion-exclusion principle, we obtain that
\begin{eqnarray*}
R_{T,p}(v) & = &
\sum_{d=1}^\infty \mu(d) M_{T,p,d}(v)
=
\sum_{\substack{d=1\\p\nmid d}}^\infty \mu(d) M_{T/d,p,1}(v)\\
& = &
\sum_{1 \le d < p} \mu(d) M_{T/d,p,1}(v)
+
\sum_{\substack{  d \ge p\\p\nmid d}} \mu(d) M_{T/d,p,1}(v)\\
& = &
\sum_{1 \le d <p } \mu(d)\frac{(T/d)^2}{p} \\
& &
\qquad +  O\(\sum_{1 \le d  < p}\left|M_{T/d,p,1}(v) -
\frac{(T/d)^2}{p}\right|
+ \sum_{ d \ge p} M_{T/d,p,1}(v)\).
\end{eqnarray*}
We see that
\begin{eqnarray*}
\sum_{1 \le d < p} \mu(d)\frac{(T/d)^2}{p}
   & = & \frac{T^2}{p} \(\zeta(2)^{-1} + O(p^{-1})\) =
\frac{T^2}{p} \( \frac{6}{\pi^2} + O(p^{-1})\),
\end{eqnarray*}
where $\zeta(s)$ is the Riemann zeta function.
Therefore
\begin{equation}
\label{eq:Deltas}
\sum_{1 \leq v \leq p-1  }
\left|R_{T,p}(v) - \frac{6}{\pi^2} \cdot \frac{T^2}{p}\right|
= O(T^2 p^{-1} + \Delta_1  + \Delta_2),
\end{equation}
where
\begin{eqnarray*}
\Delta_1
   & = & \sum_{1 \le d  < p} \sum_{1 \leq v \leq p-1 }
   \left|M_{T/d,p,1}(v) -
\frac{(T/d)^2}{p}\right|,  \\
\Delta_2
   & = & \sum_{d\ge p} \sum_{1 \leq v \leq p-1} M_{T/d,p,1}(v)  .
\end{eqnarray*}

Using the Cauchy inequality, we deduce that
\begin{equation}
\label{eq:Cauchy}
\(\sum_{1 \leq v \leq p-1 }
\left|M_{T/d,p,1}(v) - \frac{(T/d)^2}{p}\right|\)^2
\le
p
\sum_{1 \leq v \leq p-1} \left|M_{T/d,p,1}(v) - \frac{(T/d)^2}{p}\right|^2.
\end{equation}
We now recall  the bound
\begin{equation}
\label{eq:Aver Square}
\sum_{1 \leq v \leq p-1} \left|M_{W,p,1}(v) - \frac{W^2}{p}\right|^2
\le W^2 p^{o(1)},
\end{equation}
which is a special case of more general results of~\cite{Shp}
(note that the results of~\cite{Shp} apply to
the congruence $\alpha \equiv v \beta \pmod p$ where $\beta$ is not necessarily
relatively prime to $p$, but the difference of $O\(W^2/p^2\)$ for each $v$
does not affect the total error term).
We now derive from~\eqref{eq:Cauchy} and~\eqref{eq:Aver Square} that
\begin{equation}
\label{eq:Delta1}
\Delta_1 \le \sum_{1 \le d < p} \sqrt{p^{1 + o(1)}(T/d)^2}
= p^{1/2 + o(1)} T\sum_{1 \le d < p} \frac{1}{d}
= p^{1/2 + o(1)} T.
\end{equation}

The trivial bound
$$
   \sum_{1 \leq v \leq p-1} M_{W,p,1}(v)   \le W^2
$$
implies that
\begin{equation}
\label{eq:Delta2}
\Delta_2 \le \sum_{d > p} (T/d)^2 = O(T^2p^{-1}).
\end{equation}
Substituting~\eqref{eq:Delta1} and~\eqref{eq:Delta2}
in~\eqref{eq:Deltas},
we derive the desired result.

\section{Proof of Theorem~\ref{thm:Aver L-T}}

\subsection{Preliminaries}

\hspace*{0.5cm}
For a fixed $\tau \in \Q$, let $E(\tau)$ denote
the elliptic curve over $\Q$ obtained by specializing
$E(t)$ at $t = \tau$.
Let $\Delta(\tau)$ and $N(\tau)$ denote its
discriminant and conductor, respectively.
For a prime $p \nmid N(\tau)$,
let $E_p(\tau)$ denote the reduction of $E(\tau)$ modulo $p$, and let
$a_p(\tau) = p+1-\#E_p(\tau)$. Without loss of generality,
we assume that $p \geq 5$.

Now let  $\ell \neq p$ be primes such that $\ell \geq 17$ and
$j(t) $ is non-constant in $\F_p(t)$.
Let $\L = \F_p(t)$
and let $[\L]$ be its
set of places. Let
$\cB \subseteq [\L]$ be the set of places of bad reduction of
  $E(t)/\L$, which is finite
and has the property that $\deg\cB$ is bounded by a constant
independent of $p$.
Let $\L(E(t)[\ell])/\L$  be the extension of $\ell$-division points
of $E(t)$. Since $p \geq 5$,  this is a
tamely ramified Galois extension, whose Galois group we denote  $G_{\ell}$.

Since $\ell \geq 17$, we know from~\cite[Theorem~1]{CojHal} that the
geometric Galois
group of $\L(E[\ell])/\L$ is $\SL_2(\Z/\ell \Z)$.  Equivalently,
the Galois group $G_{\ell}$ is the unique subgroup of
$\GL_2(\Z/\ell \Z)$ containing $\SL_2(\Z/\ell \Z)$ and satisfying
$\det(G_{\ell}) = \langle p \rangle$.

We set
$$
G_{\ell}^{p} = \{g \in G_{\ell}: \det(g) = p\}
$$
and
$$
C^p = C \cap G_{\ell}^{p}
$$
for a finite union $C$ of conjugacy classes of $G_{\ell}$.
We have the following particular case of Murty and
Scherk~\cite[Theorem~2]{MurSch} (see also Section~3 of~\cite{CojHal}).

\begin{lemma}
\label{lem:chebotarev}
Let $U \subseteq [\L]$ be the open complement of the ramification locus $Z
\subseteq [\L]$  of $\L(E(t)[\ell])/\L$. For $v \in U(\F_p)$, let
$\Frob_v$ denote
the Frobenius at $v$ in
$\L(E(t)[\ell])/\L$. Then
$$
\#\{v \in U(\F_p): \Frob_v \subseteq C^{p}\} =
\frac{|C^p|}{\ell \left(\ell^2 -1\right)} |U(\F_p)| + O_{g, d}\left(|C^p|^{1/2}
p^{1/2}\right),
$$
where the implied $O_{g, d}$-constant depends only on the genus $g$
of $\L$ and the degree $d$ of $Z$.
\end{lemma}

  We  use this result to prove  Theorem~\ref{thm:Aver L-T}.

For Part~3 we also need the following elementary  result
  (see~\cite[Lemma~14]{CojDav}, for example).

  \begin{lemma}
\label{lem:Poly}
  Let $a, b$ be independent variables.  Let $h, w \geq 1$ be integers.
  Then there exists a polynomial
$P(X) \in \Z[X]$ such
  that
  $$
  \frac{(a^{h w} + b^{h w})^2}{ (a b)^{h w}} = P\left(\frac{(a+b)^2}{ab}\right).
  $$
  \end{lemma}

\subsection{Parts~1 and 2: Frobenius traces}

\hspace*{0.5cm}
To prove Part~1,
we follow the same lines as in the proof of~\cite[Theorem~4]{CojHal}.
In particular, for any prime $\ell$ we have that
\begin{eqnarray*}
\lefteqn{\sum_{\substack{\tau \in \cF(T)\\ \Delta(\tau) \ne 0}}
\Pi_{E(\tau)}(a, x) =
\sum_{\substack{\tau \in \cF(T)\\ \Delta(\tau) \ne 0}}\,
\sum_{\substack{p \le x, \\ p\nmid N(\tau)\\ a_{p}(\tau) =a}}1 \le
\sum_{\substack{\tau \in \cF(T)\\ \Delta(\tau) \ne 0}}\,
\sum_{\substack{p \le x, \\ p\nmid N(\tau)\\ a_{p}(\tau)
\equiv  a \pmod \ell}}1} \\
& & \qquad \qquad\le  \sum_{p\le x}\,
  \sum_{\substack{\tau \in \cF(T)\\ \Delta(\tau) \ne 0\\
N(\tau)\not \equiv
0 \pmod p \\ a_{p}(\tau) \equiv  a \pmod \ell}} 1
=\sum_{p\le x}\,
\sum_{\substack{v=1\\ \Delta(v)  N(v)\not \equiv
0 \pmod p \\a_p(v) \equiv  a \pmod \ell}}^{p-1} R_{T,p}(v).
\end{eqnarray*}

Now, applying Theorem~\ref{thm:R(T)}, we obtain
\begin{eqnarray*}
\lefteqn{
\sum_{\substack{v=1\\ \Delta(v)  N(v)\not \equiv
0 \pmod p \\a_p(v) \equiv  a \pmod \ell}}^{p-1} R_{T,p}(v) }\\
\quad & & \le
\sum_{\substack{v=1\\ \Delta(v)  N(v)\not \equiv
0 \pmod p \\a_p(v) \equiv  a \pmod \ell}}^{p-1}
\frac{6 }{\pi^2} \cdot \frac{T^2}{p}+
\sum_{\substack{v=1\\ \Delta(v)  N(v)\not \equiv
0 \pmod p \\a_p(v) \equiv  a \pmod \ell}}^{p-1}
\left|R_{T,p}(v) - \frac{6 }{\pi^2} \cdot \frac{T^2}{p}\right| \\
\quad & & \le \frac{6 }{\pi^2} \cdot \frac{T^2}{p}
\sum_{\substack{v=1\\ \Delta(v)
N(v)\not
\equiv 0 \pmod p \\a_p(v) \equiv  a \pmod \ell}}^{p-1}  1+
\sum_{v=1}^{p-1} \left|R_{T,p}(v) - \frac{6 }{\pi^2}
\cdot \frac{T^2}{p}\right| \\
\quad & & \le \frac{6 }{\pi^2} \cdot \frac{T^2}{p}
\sum_{\substack{v=1\\ \Delta(v) N(v)\not
\equiv 0 \pmod p \\a_p(v) \equiv  a \pmod \ell}}^{p-1}  1+
O\(T^2p^{-1} + T p^{1/2 + o(1)}\).
\end{eqnarray*}
Hence,
\begin{eqnarray*}
\lefteqn{\sum_{\substack{\tau \in \cF(T)\\ \Delta(\tau) \ne 0}}
\Pi_{E(\tau)}( a, x)}\\
  & & \qquad\quad \le
\frac{6T^2}{\pi}
\sum_{p\le x}\frac{1}{p}
\sum_{\substack{1\leq v \leq p-1\\
\Delta(v)  N(v)\not \equiv 0 \pmod p  \\ a_{p}(v)
\equiv  a \pmod \ell}} 1 +
O\(\sum_{p\le x}\(T^2p^{-1} + T p^{1/2 + o(1)}\)\).
\end{eqnarray*}
Using Lemma~\ref{lem:chebotarev} as  in~\cite[Theorem~2]{CojHal}
with $C$ equal to
$$
C_{\ell} = \{g \in G_{\ell}\ : \ \tr(g) = a\},
$$
we obtain that the inner sum over $v$
is  $p/\ell + O(\ell p^{1/2})$ (provided that $\ell \ge 17$).
Therefore
$$
\sum_{\substack{\tau \in \cF(T)\\ \Delta(\tau) \ne 0}}
\Pi_{E(\tau)}(a, x) \ll T^2 x\ell^{-1} + \ell T^2 x^{1/2} + T x^{3/2 + o(1)}.
$$
Finally,  by choosing $\ell$ as the smallest prime
with $\ell \ge \max\{17, x^{1/4}\}$, we conclude the proof
of Part~1 of Theorem~\ref{thm:Aver L-T}.

To prove Part 2, we remark that the condition $\tr(g) = 0$ defining $C_{\ell}$
makes sense not only in $\GL_2(\Z/\ell \Z)$, but also in $\PGL_2(\Z/\ell \Z)$.
Therefore we apply Lemma~\ref{lem:chebotarev} to the field extension
corresponding to the projection of $G_{\ell}$ in $\PGL_2(\Z/\ell \Z)$. Then
$$
 \sum_{\substack{1\leq v \leq p-1\\
\Delta(v)  N(v)\not \equiv 0 \pmod p  \\ a_{p}(v)
\equiv  0 \pmod \ell}} 1
=
p/\ell + O(\ell^{1/2} p^{1/2})
$$
(provided  that $\ell \geq 17$), and so
$$
\sum_{\substack{\tau \in \cF(T)\\ \Delta(\tau) \ne 0}}
\Pi_{E(\tau)}(0, x) \ll T^2 x\ell^{-1} + \ell^{1/2} T^2 x^{1/2} + T x^{3/2 + o(1)}.
$$
By choosing $\ell$ as the smallest prime
with $\ell \ge \max\{17, x^{1/3}\}$, we conclude the proof
of Part~2 of Theorem~\ref{thm:Aver L-T}.

\subsection{Part~3: Frobenius fields}

\hspace*{0.5cm}
As in Part~1, we have that
\begin{eqnarray*}
\lefteqn{\sum_{\substack{\tau \in \cF(T)\\
  \Delta(\tau) \ne 0
  }
  }
  \Pi_{E(\tau)}( \K, x)
=
\sum_{p \le x}
\sum_{
\substack{
  1 \leq v \leq p-1
  \\
  \Delta(v) N(v) \not \equiv 0 \pmod p
  \\
  a_p(v) \neq 0
  \\
\Q\left(\sqrt{a_p(v)^2 - 4 p}\right) = \K
  }
  }
R_{T,p}(v)}
\\
&& \qquad =
\frac{6T^2}{\pi}
\sum_{p\le x}\frac{1}{p}
   \sum_{\substack{0\leq v \leq p-1\\ \Delta(v) N(v)\not \equiv
0 \pmod p \\
a_p(v) \neq 0
\\
\Q\left(\sqrt{a_p(v)^2 - 4 p}\right) = \K
}}
1
+ O\(\sum_{p\le x}\(T^2p^{-1} + T p^{1/2 + o(1)}\)\).
\end{eqnarray*}
It remains to estimate the inner sum in the first term.

Let $0 \leq v \leq p-1$ be such that
$p \nmid \Delta(v), p \nmid N(v)$, $a_p(v) \neq 0$ and
$\Q\left(\sqrt{a_p(v)^2 - 4 p } \right) = \K$.
Let $\pi_p(v)$ be defined by
$$
X^2 - a_p(v) X + p =
\left(X - \pi_p(v)\right) \left(X - \overline{\pi_p(v)}\right).
$$
Then
$
\K
=
\Q\left(\sqrt{a_p(v)^2 - 4 p } \right)
=
\Q(\pi_p(v)),
$
and so
$p$ splits completely in $\K$. We write
$
p {\cO}_\K = \mathfrak{p} \overline{\mathfrak{p}}
$
for some conjugate prime ideals $\mathfrak{p},
\overline{\mathfrak{p}}$ of ${\cO}_\K$. In
particular, $\mathfrak{p} = (\pi_p(v))$.

Let $h$ and $w$ be the class number and the number of units of ${\cO}_\K$.
We define
$$
\pi_{p}(\K) \in {\cO}_\K
$$
by
$$
\pi_p(\K) = \alpha^w, \; \; \text{where} \; \mathfrak{p}^h = \alpha {\cO}_\K.
$$
(Note that we have two choices for $\pi_p(\K)$, and we simply make one.)
By combining the above observations,
we obtain that
\begin{equation}\label{condition1}
\pi_p(v)^{h w} = \pi_p(\K).
\end{equation}

We reinterpret~\eqref{condition1} as a Chebotarev
condition in some extension of $\F_p(t)$ (note that here $p$ and $\K$
are fixed,
and $v$ is a specialization of $t$). To do this, let us choose a rational prime
$\ell \geq 17$, $\ell
\neq p$, and consider
the Galois extension $\L(E[\ell])/\L$.
 From classical theory we know that the  Frobenius at $v$ in this extension,
  viewed
  as an element of $\GL_2(\Z/\ell \Z)$, has the property that its trace
$\tr \Frob_v $ satisfies
  $$
  \tr \Frob_v \equiv \pi_p(v) + \overline{\pi_p(v)} \pmod \ell.
  $$
  Thus $\pi_p(v)$ has a ``Chebotarev interpretation'' in this extension.

  Now we combine~\eqref{condition1} with  Lemma~\ref{lem:Poly}, getting
  $$
\frac{\left(\pi_p(\K) + \overline{\pi_p(\K)}\right)^2}{
\pi_p(\K) \overline{\pi_p(\K)}}
=
\frac{\left(\pi_p(v)^{h w} + \overline{\pi_p(v)}^{h w}\right)^2}{
\pi_p(v)^{h w} \overline{\pi_p(v)}^{h w}
}
=
P\left(
\frac{(\pi_p(v) + \overline{\pi_p(v)})^2}{\pi_p(v) \overline{\pi_p(v)}}
\right).
$$
Let us define
$$
C_{\ell} =
\left\{
g \in G_{\ell}\ : \
P\(\frac{\Tr(g)^2}{\det g}\right)
=
\frac{\left(\pi_p(\K) + \overline{\pi_p(\K)}\right)^2}{
\pi_p(\K) \overline{\pi_p(\K)}}
\right\}.
$$
where $\Tr(g)$ and $\det g$ denote the trace and determinant of $g$, respectively.
Then
\begin{equation}
\label{condition2}
  \sum_{\substack{0\leq v \leq p-1\\ \Delta(v) N(v)\not \equiv
0 \pmod p \\
a_p(v) \neq 0
\\
\Q\left(\sqrt{a_p(v)^2 - 4 p}\right) = \K
}}
1
  \leq
\#\left\{
1 \leq v \leq p-1 \ : \
   p \nmid \Delta(v) N(v), \
\Frob_v \subseteq C_{\ell}
\right\}.
\end{equation}

To estimate~\eqref{condition2} we can now invoke Lemma~\ref{lem:chebotarev}.
Again, as in the proof of Part 2, we remark that the condition defining $C_{\ell}$
makes sense not only in $\GL_2(\Z/\ell \Z)$, but also in $\PGL_2(\Z/\ell \Z)$.
Thus we apply Lemma~\ref{lem:chebotarev} to the field extension corresponding
to the projection of $G_{\ell}$ in $\PGL_2(\Z/\ell \Z)$.
It is an easy calculation to show that $\# C_{\ell}^p =O(\ell^2)$ and
 $\# \overline{C}_{\ell}^p = O(\ell)$, where $\overline{C}_{\ell}^p$ is the
 union of conjugacy classes in $\PGL_2(\Z/\ell \Z)$ 
 of the elements in $C_{\ell}^{p}$.
After putting everything together and continuing as in Part~2,
we conclude the proof.


\end{document}